\documentclass[centertags,12pt]{amsart}
\usepackage{latexsym}
\usepackage{amssymb}
\usepackage{graphics}
\numberwithin{equation}{section}
\makeatletter
\renewcommand{\subsection}{\@startsection
{subsection}{2}{0mm}{\baselineskip}{-0.25cm}
{\normalfont\normalsize\rm}}
\makeatother
\newtheorem{theorem}{Theorem}[section]
\newtheorem{proposition}[theorem]{Proposition}

\newtheorem{lemma}[theorem]{Lemma}
{\theoremstyle{definition}

\newtheorem{example}[theorem]{Example}

}
\theoremstyle{remark}
\newtheorem{remark}[theorem]{Remark}
\newtheorem*{claim*}{Claim}
\def\1{\mathbf 1}

\def\N{\mathbf N}

\def\cC{\mathcal C}

\def\cF{\mathcal F}

\def\cL{\mathcal L}
\def\cP{\mathcal P}
\def\cX{\mathcal X}

\def\supp{{\rm supp}}
\def\div{{\rm div}}
\def\Div{{\rm Div}}
\def\dim{{\rm dim}}
\def\deg{{\rm deg}}

\def\ltilde{\tilde{\ell}}
\def\fq{{\mathbf F_q}}
\def\fqs{{\mathbf F_{q^2}}}
\def\bx{\mathbf x}
\begin{document}
\author[C.~Munuera]{Carlos Munuera}
\author[F.~Torres]{Fernando Torres}
\thanks{The authors were partially supported respectively by the Grants
VA020-02 (``Junta de Castilla y Le\'on"), Proc. 300681/97-6 (CNPq-Brazil)  
and SB2000-0225 (``Secretaria de Estado de Educaci\'on y Universidades del
Ministerio de Educaci\'on, Cultura y Deportes de Espa\~na")}
\title[Bounding the trellis state complexity of AG codes]{Bounding the
trellis state complexity of\\ algebraic geometric codes}
\address{Dept. of Applied Mathematics, University of Valladolid (ETS
Arquitectura), Avda. Salamanca SN, 47014 Valladolid, Castilla, Spain}
\email{cmunuera@modulor.arq.uva.es}
\address{IMECC-UNICAMP, Cx. P. 6065, Campinas, 13083-970, SP - Brazil}
\email{ftorres@ime.unicamp.br} 
     
    \begin{abstract} Let $\cC$ be an algebraic geometric code of dimension
$k$ and length $n$ constructed on a curve $\cX$ over $\fq$. Let $s(\cC)$
be the state complexity of $\cC$ and set $w(\cC):=\min\{k,n-k\}$, the Wolf
upper bound on $s(\cC)$. We introduce a numerical function $R$ that
depends on the gonality sequence of $\cX$ and show that $s(\cC)\geq
w(\cC)-R(2g-2)$, where $g$ is the genus of $\cX$. As a matter of fact,
$R(2g-2)\leq g-(\gamma_2-2)$ with $\gamma_2$ being the gonality over $\fq$
of $\cX$, and thus in particular we have that $s(\cC)\geq
w(\cC)-g+\gamma_2-2$.
    \end{abstract}
    
\maketitle
   \section{Introduction}\label{s1}
A {\em trellis of depth $n$} is an edge-labeled directed graph $T=(V,E)$
with vertex set $V$ and edge set $E$ satisfying the following properties:
   \begin{itemize}
\item $V$ is the union of $(n+1)$ disjoint subsets $V_0, \ldots,V_n$;
\item every edge in $E$ that begins at $V_i$ ends at $V_{i+1}$;
\item every vertex in $V$ belongs to at least one path from a vertex in
$V_0$ to a vertex in $V_n$.
   \end{itemize} 
In this paper we only consider trellises with $V_0$ and $V_n$ having just
one element. To each path from $V_0$ to $V_n$ one can associate an ordered
$n$-tuple over a label alphabet, say $\fq$ the finite field with $q$
elements. Thus the set of all such $n$-tuples defines a block code $\cC_T$
of length $n$ over $\fq$. Conversely, given a block code $C\subseteq
\mathbf F^n_q$ we say that a trellis $T$ {\em represents} $C$ if
$\cC_T=\cC$. There might exist more than one non-isomorphic
trellis representing the same code. The use of trellises in coding
theory started with applications to convolutional codes. Then they were
employed with block codes mainly for the purpose of soft-decision decoding
with the Viterbi algorithm. History and the state of the art of
application of trellises to coding theory can be seen in Forney's paper
\cite{forney} and Vardy's survey \cite{vardy}.
A way to measure the complexity of a trellis $T$ that represents a
code $\cC\subseteq \mathbf F^n_q$ is by means of {\em the state of
complexity} of $T$ denoted by $s_T(\cC)$ and defined by
   $$
s_T(\cC):=\max\{s_0(T),s_1(T),\ldots,s_n(T)\}\, , 
   $$
where $s_i(T):={\rm log}_q|V_i|$ with $V_0, \ldots,V_n$ being the
underlying partition of the vertex set of $T$. If the code $\cC$ is
linear, and once the order of coordinates of $\cC$ is
fixed, there exists an unique (up to a graph isomorphism) trellis
$T_\cC$ such that for each $i=0,1,\ldots,n$, and any trellis $T$ that
represents $\cC$ it holds that $s_i(T_\cC)\leq s_i(T)$. The trellis
$T_\cC$ is called {\em the minimal trellis} of $\cC$. Then the sequence 
$(s_i(\cC):i=0,1,\ldots,n)$, with $s_i(\cC):=s_i(T_\cC)$, is called the
{\em state complexity profile} of $\cC$ and the
number $s(\cC):= s_{T_\cC}(\cC)$ is the {\em state complexity} of
$\cC$. Forney \cite{forney} (see
also \cite[Ex. 5.1]{vardy}) noticed that the state complexity of $\cC$
may vary when changing the order of coordinates. We shall say that two 
codes are {\em equivalent} if one of them can be
obtained from the other by permuting coordinates and we denote by
$[\cC]$ the set of codes which are equivalent to $\cC$. We thus 
are lead to consider {\em the absolute state complexity} of $\cC$,
namely
   $$
s[\cC]:=\min\{s(\cC'): \cC' \in [\cC]\}\, .
   $$
We mention that the role of the state complexity of a 
linear code is comparable to the role that plays its length, its dimension
and its minimum distance; cf. Muder \cite{muder}, Forney \cite{forney}. 
In general, there are several bounds on $s(\cC)$ available in the
literature, see e.g. \cite[Sects. 5.2,5.3]{vardy}. Here we just mention
the Wolf bound, as it was first noticed by him in \cite{wolf}, namely
    $$
s(\cC)\leq w(\cC):=\min\{k,n-k\}\, ,
    $$
where $k$ is the dimension of $\cC$. The study of the state complexity of
some classical codes, such as BCH, RS, and RM codes, has been carried out
by several authors; see \cite{berger-beery1}, \cite{berger-beery2}, 
\cite{blackmore-norton1}, \cite{kasami-all}, \cite{vardy-beery}. The case
of algebraic geometric codes (or simply, AG codes)
was treated by Shany and Be'ery \cite{shany-beery}, Blackmore and Norton
\cite{blackmore-norton2}, \cite{blackmore-norton3}, and by Munuera and
Torres \cite{munuera-torres}. If $\cC=\cC(\cX,D,G)$ is an AG code, from these works it
follows that $s(\cC)=w(\cC)$, provided that either
$\deg(G)<\lfloor\deg(D)/2\rfloor$, or $\deg(G)>\lceil\deg(D)/2\rceil+2g-2$,
with $g$ being the genus of the underlying curve $\cX$. Otherwise we have
the so-called {\em Clifford bound}, namely $s(\cC)\geq
\lceil\deg(D)/2\rceil-g-1$. The main result in \cite{munuera-torres} is a
Goppa-like bound on $s(\cC)$,  $s(\cC)\geq w(\cC)-(g-\alpha)$ with
$\alpha$ being the abundance of the code. The particular case of Hermitian
codes have been treated in 
\cite{blackmore-norton2}, \cite{blackmore-norton3} 
and \cite{shany-beery}. The purpose of this paper is to investigate further lower bounds on
$s(\cC)$ for $\cC=\cC(\cX,D,G)$ an AG code. Our approach depends heavily on the gonality
sequence of the curve $\cX$ used to construct the code and as a matter of
facts our results subsume  the previous aforementioned lower bounds on
$s(\cC)$. We introduce a numerical function $R(N)$ which gives rise to 
our first main result Theorem \ref{thm3.1} below establishing that
$s(\cC)\geq w(\cC)-R(2\deg(G)-\deg(D))\geq w(\cC)-R(2g-2)$. Essentially this bound is the same as the one introduced by Blackmore and Norton in \cite{blackmore-norton3}, where it is called {\em second gonality bound}. In \cite{blackmore-norton2} this bound was computed by the same authors for the case of Hermitian codes and used to give a very good bound on $s(\cC)$ for these codes.  By using specific properties of $R$ (cf. Lemma \ref{lemma3.4}), here we can compute
$R(2g-2)$ (see Proposition \ref{prop3.3}) and in particular we
find that $R(2g-2)\leq g-(\gamma_2-2)$, where $\gamma_2$ is the gonality
over $\fq$ of the curve $\cX$. In this way we obtain our second main
result, Theorem \ref{thm3.2} below, which asserts that $s(\cC)\geq
w(\cC)-g+\gamma_2-2$. Furthermore we can compute $R$ for all plane (nonsingular) curves, and thus we can extend some results of \cite{blackmore-norton2} to all codes coming from these curves.
Section \ref{s2} contains preliminary results on the state complexity of linear
codes and gonality sequence of curves. The novelty here is a
``symmetric-like property" (Proposition \ref{prop2.2}) of the gonality
sequence of a curve which has been noticed first in
\cite{carvalho-pellikaan-torres}. Finally, in Section 4 we state a new  property of Self-orthogonal codes which was first noticed by Blackmore and Norton in \cite{blackmore-norton2} for the case of Hermitian codes (see Proposition 4.4 here).
   \section{Preliminaries}\label{s2}
In this section we shall point out some results on the state
complexity of linear codes as well as some basic properties of the
gonality sequence of curves which play a role in the present work.
   \subsection{On the state complexity of linear codes}\label{s2.1}
For a $[n,k,d]$ linear code $\cC$ over the finite field $\fq$, its minimal trellis $T=T_\cC$ can be constructed in several ways;
see \cite[Sect. 4]{vardy}. For our purpose the relevant construction is the one given by
Forney. He shows that the sub-sets $V_0, \ldots,V_n$ of the
underlying partition of the vertex set of $T$ are given by
$V_i=\cC/\cP_i\oplus \cF_i$, with $\cP_i=\cP_i(\cC)$ and
$\cF_i=\cF_i(\cC)$ being respectively the $i$-th past and the $i$-th
future subcodes of $\cC$; namely, $\cP_0=\cF_n=0$, $\cP_n=\cF_0=\cC$,
and for $i=1,\ldots,n-1,$
   \begin{align*}
\cP_i & =  \{(c_1,\ldots,c_i): (c_1,\ldots,c_i,0,\ldots,0)\in\cC\}\,
,\\         
\cF_i & =  \{(c_{i+1},\ldots,c_n): (0,\ldots,0,c_{i+1},\ldots,c_n)
\in\cC\}\, .                                   
    \end{align*}
It follows that $|V_i|$ is a power of $q$ so that
   $$
s_i(\cC)=k-\Delta_i\, ,
   $$
where
   $$
\Delta_i=\Delta_i(\cC):=\dim_\fq(\cP_i)+\dim_\fq(\cF_i)\, ,
\qquad i=0,1,\ldots,n\, .
   $$ 
 
    \begin{lemma}{\rm (\cite[Thm. 4.20]{vardy})}\label{lemma2.1} The
state complexity of a linear code and that of its dual are
identical.
   \end{lemma}
In particular, from this Lemma and the Forney's construction, the Wolf bound,
$s(\cC)\leq w(\cC)=\min\{k,n-k\}$, follows. 
Therefore to study the state complexity of $\cC$ we can restrict
ourselves to the case $2k\leq n$. Now by definition
   $$
s(\cC)=k-\Delta\, ,
   $$
where $\Delta=\Delta(\cC):=\min\{\Delta_0,\Delta_1,\ldots,\Delta_n\}$. 
We set $\Delta[\cC]:=\max\{\Delta(\cC'): \cC'\in [\cC]\}$.
    \begin{lemma}\label{lemma2.2} With the above notation$,$ the following
holds$:$
   \begin{enumerate}
\item[\rm(1)] $\cP_i=0$ for $i=0,\ldots, d-1$; in particular     $\min\{\Delta_0,\ldots,\Delta_{d-1}\}=\Delta_{d-1}$.
\item[\rm(2)] $\cF_i=0$ for $i=n-d+1, \ldots,n$; in particular 
$\min\{\Delta_{n-d+1},\ldots,\Delta_n\}=\Delta_{n-d+1}$.
   \end{enumerate}
   \end{lemma}    
   \begin{proof} For $i\in\{1,\ldots,d-1\}$ (resp. $i\in
\{n-d+1,\ldots,n-1\}$), the length of $\cP_i$ (resp. $\cF_i$) is smaller
than the weight of any nonzero codeword in $\cC$. Thus $\cP_i=0$ (resp.
$\cF_i=0$) and the result follows taking into account the fact that
$\dim_\fq(\cF_0)>\dim_\fq(\cF_1)>\ldots$ (resp.
$\dim_\fq(\cP_n)>\dim_\fq(\cP_{n-1})>\ldots$).
   \end{proof}
   \begin{proposition}\label{prop2.1} For a linear code $\cC$, we have  
   $$
s(\cC)=\begin{cases}
   k  & \text{whenever $2d\geq n+2$,}\\
   k-\min\{\Delta_{d-1},\ldots,\Delta_{n-d+1}\} & 
                                         \text{otherwise.}
   \end{cases}
   $$
   \end{proposition}
   \begin{proof} The result follows from Lemma \ref{lemma2.2} by taking
into consideration that there exists an integer $i$ with $n-d+1\leq i\leq
d-1$ whenever $2d\geq n+2$.
   \end{proof}
   \subsection{On the gonality sequence of curves}\label{s2.2}
Let $\cX$ be a (projective, geometrically irreducible, non-singular
algebraic) curve defined over the finite field $\fq$. Let $i$
be a positive integer. The {\em $i$-th gonality number over $\fq$} of
$\cX$ is defined by
   $$
\gamma_i=\gamma_i(\cX,\fq):=\min\{\deg(A): \text{$A\in\Div(\cX,\fq)$ with  
                               $\ell(A)\geq i$}\}\, .
   $$ 
As usual $\Div(\cX,\fq)$ denotes the set of $\fq$-divisors of $\cX$
and for a divisor $F$, $\cL(F)$ stands for the $\fq$-vector space of
$\fq$-rational functions $f$ on $\cX$ such that $f=0$ or $F+\div(f)\succeq
0$. We set $\ell(F):=\dim_\fq\cL(F)$. Standard references on algebraic
geometry and algebraic function fields are the books by Hartshorne
\cite{hartshorne} and Stichtenoth \cite{sti} respectively.
The sequence $GS(\cX)=GS(\cX,\fq):=(\gamma_i:i\in\N)$ is called {\em the 
gonality sequence} of $\cX$ over $\fq$. Notice that $\gamma_1=0$ and that
$\gamma_2$ is the usual gonality of $\cX$ over $\fq$.
   \begin{remark}\label{rem2.1} Pellikaan \cite{pellikaan0} noticed the
relevance of the gonality of the underlying curve in the study of AG codes. 
The invariant $GS(\cX)$ was introduced by Yang, Kumar
and Stichtenoth in \cite{yang-kumar-sti} in connection with lower bounds
on the generalized Hamming weight hierarchy of AG codes; these lower
bounds were generalized by Munuera in \cite{munuera}.
    \end{remark}
Some properties of $GS(\cX)$ are stated below. Let $g$ be the genus of $\cX$. 
   \begin{lemma}{\rm (\cite[Prop. 11]{yang-kumar-sti})}\label{lemma2.3}
Suppose that $\cX(\fq)\neq \emptyset.$ Then$:$
   \begin{enumerate}
\item[\rm(1)] The sequence $GS(\cX)$ is strictly increasing$;$
\item[\rm(2)] $2i-2\leq \gamma_i\leq g+i-2,$ for $i=2,\ldots,g;$
\item[\rm(3)] $\gamma_g=2g-2;$
\item[\rm(4)] $\gamma_i=g+i-1,$ for $i\geq g+1.$
   \end{enumerate}
   \end{lemma}
In general, it is quite difficult to compute the sequence $GS(\cX)$;  nevertheless it is available in the following cases.
   \begin{lemma}\label{lemma2.4} Suppose that $\cX(\fq)\neq \emptyset.$
   \begin{enumerate}
\item[\rm(1)] If $\cX$ is a hyperelliptic curve$,$ then $\gamma_i=2i-2$ for $i=1,\ldots,g-1;$
\item[\rm(2)]{\rm (\cite[Cor. 2.4]{pellikaan})} If $\cX$ is a (non-singular) plane curve of degree $r+1,$ then $GS(\cX)$ is the strictly increasing sequence obtained from the semigroup generated by $r$ and $r+1.$
   \end{enumerate}
   \end{lemma}
   \begin{remark}\label{rem2.2} By Clifford's theorem, Item (1) in the
above result can be improved by observing that the curve $\cX$ is
hyperelliptic if and only if $\gamma_i=2i-2$ for some $i\in\{2,\ldots,g-1\}$.
    \end{remark}
We shall need the following result which was originality noticed in
\cite{carvalho-pellikaan-torres}; we include the proof for the sake of
completeness.
    \begin{proposition}\label{prop2.2} Assume $\cX(\fq)\neq\emptyset.$
Let $a$ be an integer with $0\leq a\leq 2g-1.$ Then $a\in GS(\cX)$ if and
only if $2g-1-a\not\in GS(\cX).$
    \end{proposition}
    \begin{proof} By Lemma \ref{lemma2.3} in the interval $[0,2g-1]$ there
are precisely $g$ gonality numbers of $\cX$. Thus it is enough to show
that $2g-1-\gamma_i\neq \gamma_{j}$ for any $i,j=1,\ldots, g$. Let
$A\in\Div(\cX,\fq)$ such that $\deg(A)=\gamma_i$ and $\ell(A)\geq i$. Let
$W$ be a canonical divisor on $\cX$. By the Riemann-Roch theorem,
$\ell(W-A)\geq i+g-\gamma_i-1$.
Suppose that $j\leq i+g-\gamma_i-1$. Then $\gamma_{j}\leq
\deg(W-A)=2g-2-\gamma_i$ and thus $2g-1-\gamma_i\neq \gamma_{j}$. Now let
$j\geq i+g-\gamma_i$ and suppose by means of contradiction that
$2g-1-\gamma_i=\gamma_{j}$. Let $B\in \Div(\cX,\fq)$ such that
$\deg(B)=\gamma_{j}$ and $\ell(B)\geq j$. As above we have
$\ell(W-B)\geq j+g-\gamma_{j}-1$ so that $\ell(W-B)\geq j+\gamma_i-g\geq
i$. This is not possible since $\deg(W-B)=2g-2-\gamma_{j}=\gamma_i-1$.
    \end{proof}
    \section{A lower bound on the absolute state complexity of an AG
code}\label{s3}
The goal of this section is to state and prove new lower bounds on the
absolute state complexity of an AG code (see Theorems \ref{thm3.1} and
\ref{thm3.2} below) which will be related to the gonality sequence of the
underlying curve. Standard references on AG codes are the survey \cite{hoholdt-lint-pellikaan} by H\o holdt, van Lint and Pellikaan, and Stichtenoth's book \cite{sti}.
Let $\cX$ be a (projective, geometrically irreducible, non-singular
algebraic) algebraic curve of genus $g$ defined over the field $\fq$. Let
$D$ and $G$ be two $\fq$-rational divisors on $\cX$ with
$D=P_1+\ldots+P_n$ being the sum of $n$ pairwise different $\fq$-rational
points on $\cX$ such that $P_i\not\in\supp (G)$. The AG code
$\cC=\cC(\cX,D,G)$
is the image in $\mathbf F^n_q$ of the $\fq$-linear map
   $$
ev: \cL(G)\to \mathbf F^n_q\, ,\qquad f\mapsto (f(P_1),\ldots,f(P_n)).
   $$
Notice that the kernel of $ev$ is $\cL(G-D)$. The number $\ell(G-D)$ is called {\em the 
abundance} of $\cC$ and $\cC$ is called {\em non-abundant} provided that
$\ell(G-D)=0$.
Let $k$ and $d$ be respectively the dimension and the minimum distance of
$\cC$. We have the so-called {\em Goppa estimates} on the parameters $k$
and $d$, namely
    $$
k=\ell(G)-\ell(G-D)\, ,\qquad\text{and}\qquad d \geq n-\deg(G)\, ,
   $$
and thus they can be handled by means of the Riemann-Roch theorem. 
    \begin{lemma}\label{lemma3.1} If $2k\leq n$ and $n>2g,$ then the AG
code $\cC=\cC(\cX,D,G)$ is non-abundant and $2\deg(G)-n\leq 2g-2.$
   \end{lemma}
   \begin{proof} Suppose that $\ell(G-D)\geq 1$. Then the divisor $G-D$
must be special; otherwise $k=\ell(G)-(\deg(G-D)+1-g)\geq
(\deg(G)+1-g)-(\deg(G-D)+1-g)=n$ which is a contradiction. By Clifford's
theorem, $\deg(G-D)\leq (\deg(G)-n)/2+1$ and hence
   $$
k=\ell(G)-\ell(G-D)\geq (\deg(G)+1-g)-(\deg(G)-n)/2-1=
(\deg(G)+n-2g)/2\, .
   $$ 
{}From the hypothesis $2k\leq n$ we conclude that $2g\geq \deg(G)$. On
the other hand, $\deg(G-D)\geq 0$ as $\ell(G-D)\geq 1$, and thus $2g\geq
\deg(G)\geq n$; a contradiction. The second statement follows from the
fact that $n/2\geq k=\ell(G)\geq \deg(G)+1-g$.
   \end{proof}
   \begin{remark}\label{rem3.0} Examples of curves of genus $g$ over $\fq$
having more than $2g$ $\fq$-rational points are the maximal
curves over $\fq$, namely those whose number of $\fq$-rational points
attains the Hasse-Weil upper bound. Numerical examples are the Fermat
curves over $\fqs$ of degree $r+1$, a divisor of $q+1$. Other examples are
the Hurwitz curves $X^hY+Y^hZ+Z^hX=0$ over $\fqs$, where
$(h^2-h+1)$ is a divisor of $q+1$; see e.g. \cite{akt}.
   \end{remark}
Let $GS(\cX)=(\gamma_i:i\in\N)$ be the gonality sequence over $\fq$ of
$\cX$. It will be convenient for us to consider $GS(\cX)$ as a subset of
$\N':=\{-1\}\cup\N_0$. An element in $\N'\setminus GS(\cX)$ will be called
a {\em gap} of $\cX$. By Lemma \ref{lemma2.3} there are $g+1$ gaps and the
biggest one is $2g-1$.
Let $\tilde\ell=\tilde\ell_\cX:\N'\to\N_0$ be the numerical function
defined by $\tilde\ell(-1):=0$ and
   $$
\tilde\ell(a):=\max\{i\in\N: \gamma_i\leq a\}\, ,\qquad a\in\N_0\, .
   $$
By Lemma \ref{lemma2.3} the function $\ltilde$ becomes an increasing step
function such that $\ltilde(2g-2)=g$ and $\ltilde(2g-1+i)=g+i$ for $i\geq
0$. Moreover, $\ltilde(a+1)\leq \ltilde(a)+1$ and equality holds if and
only if $a+1\in GS(\cX)$.
   \begin{lemma}\label{lemma3.2} For $F$ a rational divisor on $\cX$ with
$\deg(F)\geq -1,$ $\ell(F)\leq \ltilde(\deg(F)).$
   \end{lemma}
   \begin{proof} If $\deg(F)=-1$, then $\ell(F)=0=\ltilde(-1)$. Let
$\deg(F)\geq 0$ and let $i\in\N_0$ be such that $\gamma_i\leq
\deg(F)<\gamma_{i+1}$ so that $\ltilde(\deg(F))=i$. Thus by definition of
$\gamma_{i+1}$ we must have $\ell(F)\leq i$ and the result follows.
   \end{proof}
Next we let $R=R_\cX:\N'\cap[-1,2g-2]\to \N$ be the numerical function
defined by
   $$
R(N):=\min\{\ltilde(a)+\ltilde(b):\text{$a,b\in \N'$ with $a+b=N$}\}\, .
   $$
Now we can state the first main result of this section.
    \begin{theorem}\label{thm3.1} Let $\cC=\cC(\cX,D,G)$ be an AG code
such that $2k\leq n$ and $n>2g,$ where $g$ is the genus of $\cX.$ Set
$m:=\deg(G).$ Then
   $$
\Delta[\cC]\leq R(2m-n)
   $$
and hence  $s[\cC]\geq w(\cC)-R(2m-n)$.
   \end{theorem}
   \begin{proof} Since the function $R$ depends only on the underlying
curve $\cX$, it is enough to show that $\Delta(\cC)\leq R(2m-n)$. In
addition, $w(\cC)=\min\{k,n-k\}=k$ and by Proposition \ref{prop2.1} we can
assume that $2d<n+2$ so that $2m-n\geq -1$ by the Goppa estimative on $d$.
Let us recall that the $i$-th past and the $i$-th future subcodes 
$\cP_i$ and $\cF_i$ respectively of the AG code $\cC$ in Forney's
construction are also AG codes and are given by (see
\cite{blackmore-norton1})     
   \begin{align*}
\cP_i & =  \cC(\cX,D-P_{i+1}-\ldots-P_n,G-P_{i+1}-\ldots-P_n)\, ,\\         
\cF_i & =  \cC(\cX,D-P_1-\ldots-P_i,G-P_1-\ldots-P_i)\, .
   \end{align*}
Now by Lemma \ref{lemma3.1} the code $\cC$ is non-abundant and hence 
the $i$-th element $s_i=s_i(\cC)$ in the state complexity profile of
$\cC$ is given by
   $$
s_i=k-\Delta_i\, ,
   $$
where $\Delta_i=\ell(G-P_1-\ldots-P_i)+\ell(G-P_{i+1}-\ldots-P_n)$. Thus,
according to Proposition \ref{prop2.1}, and since $d\geq n-m$, we have
$s(\cC)=w(\cC)-\Delta(\cC)$ where 
   $$
\Delta(\cC)=  \min\{\Delta_{d-1},\ldots,\Delta_{n-d+1}\}
           =  \min\{\Delta_{n-m-1}, \ldots,\Delta_{m+1}\}\, .
   $$
Let $i$ be an integer with $n-m-1\leq i\leq m+1$ so that
$\deg(G-P_1-\ldots-P_i)\geq -1$ and $\deg(G-P_{i+1}-\ldots-P_n)\geq -1$;
then by Lemma \ref{lemma3.2},
    $$
\Delta_i\leq \ltilde(\deg(G-P_1-\ldots-P_i))+
\ltilde(\deg(G-P{i+1}-\ldots-P_n))\, .
   $$ 
Now as $\deg(G-P_1-\ldots-P_i)+\deg(G-P_{i+1}-\ldots-P_n)=2m-n$
which is at most $2g-2$ by Lemma \ref{lemma3.1}, the result follows.
    \end{proof} 
In order to apply the above result we need to know the behavior of the
function $R$. This study is done in the rest of this section. In
particular, we shall compute $R(2g-2)$ whenever $g>0$ and also explicitely
describe $R$ for the case of plane curves.
   \begin{lemma}\label{lemma3.4} Let $N\in \N'\cap[-1,2g-2].$
   \begin{enumerate}
\item[\rm(1)] $R$ is an increasing function such that $R(N)<R(N+1)$
implies $R(N+1)=R(N)+1;$
\item[\rm(2)] If $N<\gamma_i-1,$ then $1\leq R(N)\leq i-1;$
\item[\rm(3)] $R(N)\leq \lfloor (N+1)/2\rfloor +1;$
\item[\rm(4)] There exists a gap $a=a(N)$ of $\cX$ with $a\leq N/2$ such
that $R(N)=\ltilde(a)+\ltilde(N-a).$
    \end{enumerate}
    \end{lemma}
    \begin{proof} From the definition of $R$ it is clear that $R(N)\geq 1$
and that $R(-1)=1$.
(1) Let $R(N+1)=\ltilde(a)+\ltilde(b)$ with $a+b=N+1$ and $a\leq b$. From
$a+(b-1)=N$ we have $R(N)\leq \ltilde(a)+\ltilde(b-1)\leq
\ltilde(a)+\ltilde(b)=R(N+1)$ since $\ltilde$ is an increasing function.
Now suppose that $R(N)<R(N+1)$ and let $R(N)=\ltilde(a')+\ltilde(b')$ with
$a'+b'=N$. Then from $(a'+1)+b'=N+1$, $R(N+1)\leq
\ltilde(a'+1)+\ltilde(b')$ and thus $\ltilde(a'+1)>\ltilde(a')$. Therefore
$R(N+1)=R(N)+1$ since $\ltilde(a'+1)=\ltilde(a')+1$.
(2) From $N=-1+(N+1)$ it follows that $R(N)\leq \ltilde(N+1)$; the latter
number is at most $i-1$ by hypothesis and (2) follows.
(3) There exists $i\in\{1,\ldots,g\}$ such that $\gamma_i\leq
N+1<\gamma_{i+1}$. Then by (2), $R(N)\leq i$, and the latter number is at
most $(N+3)/2$ by Lemma \ref{lemma3.2}(2).
(4) Let $R(N)=\ltilde(a)+\ltilde(b)$ with $a\leq b=N-a$ and suppose that
$a\in GS(\cX)$. We have $\ltilde(a-1)=\ltilde(a)-1$, and $\ltilde(b+1)\leq
\ltilde(b)+1$. Then $\ltilde(a-1)+\ltilde(b+1)\leq
\ltilde(a)+\ltilde(b)=R(N)$ and thus $R(N)=\ltilde(a-1)+\ltilde(b+1)$. If
$a-1$ is a gap of $\cX$, then we are done; otherwise we repeat the above
argument.
    \end{proof}
    \begin{remark}\label{rem3.1} From Lemma \ref{lemma3.4}(1)(3), we have
that $R(N)\leq R(2g-2)\leq g$ whenever $N\in\N'\cap[-1,2g-2]$. Then
Theorem \ref{thm3.1} yields the main result in \cite{munuera-torres},
namely $s[\cC]\geq w(\cC)-g$, provided that $2k\leq n$ and $n>2g$. We are
going to improve this result via Proposition \ref{prop3.3} and 
Theorem \ref{thm3.2} below.
    \end{remark}
    \begin{lemma}\label{lemma3.3} Let $i\in \N', N\in\N'\cap[-1,2g-2]$ and
$r\in\N$ with $i+r\leq N+1.$ If $A=\{i,i+1,\dots,i+r\}\subseteq \N'$ is a
set of $r+1$ consecutive integers such that $i+1,\dots,i+r$ are gaps of
$\cX$, then   
   $$
\min\{\ltilde(a)+\ltilde(N-a): a\in A \}=\ltilde(i+r)+\ltilde(N-i-r)\, ..
   $$ 
   \end{lemma}
   \begin{proof} Let $a=i+j$ with $1\leq j\leq r$. We have that
$\ltilde(a)=\ltilde(i)$ since by hypothesis $a$ is a gap of $\cX$. Then 
$\ltilde(a)+\ltilde(N-a)$ is minimum when $\ltilde(N-a)$ is; i.e., when
$a$ is the largest element in $A$ as $\ltilde$ is an
increasing function.
   \end{proof}
   
 \begin{proposition}\label{prop3.1} For $N\in\N'\cap[-1,2g-2],$ 
   \begin{align*}
R(N)= & \min\{\ltilde(a)+\ltilde(N-a): \text{$-1\leq a\leq N/2, a=\lfloor
N/2\rfloor$ or} \\
   {}   & \text{$a\in \N'\setminus GS(\cX)$ with $a+1\in GS(\cX)$}\}\, .
   \end{align*}
   \end{proposition}
   \begin{proof} By Lemma \ref{lemma3.4}(4),
$R(N)=\ltilde(a)+\ltilde(N-a)$
for some gap $a$ of $\cX$ such that $a\leq N/2$. Suppose that $a<\lfloor
N/2\rfloor$. If each integer $a'$ with $a<a'\leq \lfloor N/2\rfloor$ is a
gap of $\cX$, then from Lemma \ref{lemma3.3} $R(N)=\ltilde(\lfloor
N/2\rfloor)+\ltilde(\lceil N/2\rceil)$; otherwise, by Lemma \ref{lemma3.3}
again, we can assume $a+1\in G(\cX)$ and the result follows.
   \end{proof}
   \begin{example}\label{ex3.1} Let $N\in\N'\cap[-1,2g-2]$ such that
$\lceil N/2\rceil<\gamma_2$. Then from the above result we have that
$R(N)=1$ whenever $N+1<\gamma_2$, and $R(N)=2$ otherwise.
   \end{example}
Next we show that the upper bound for $R(N)$ in Lemma \ref{lemma3.4}(3)
is the best possible.
   \begin{proposition}\label{prop3.2} If $\cX$ is a hyperelliptic curve$,$
then
   $$
R(N)=\lfloor\frac{N+1}{2}\rfloor +1\, ,\qquad N\in\N'\cap[-1,2g-2]\, .
  $$ 
Conversely$,$ suppose the above formula holds true for some
$N\in\N'\cap[1,2g-4].$ 
   \begin{enumerate}
\item[\rm(1)] If $N$ is odd$,$ then $\cX$ is hyperelliptic$;$
\item[\rm(2)] If $N$ is even$,$ then either $\cX$ is hyperelliptic$,$ or
$\gamma_i=2i-1$ for some $i\in\{2,\ldots,g-1\}.$
   \end{enumerate}
   \end{proposition}
   \begin{proof} Let $\cX$ be hyperelliptic and $a\in 
\N'\cap[-1,2g-1]$. By Lemma \ref{lemma2.4}(1), $a$ is a gap of
$\cX$ if
and only if $a$ is odd; moreover,  $\ltilde(a)=\lfloor (a+2)/2\rfloor$.
Now let $R(N)=\ltilde(a')+\ltilde(N-a')$ with $a'$ a gap of $\cX$ (cf.
Lemma \ref{lemma3.4}(4)). Then 
   $$
R(N)=\frac{a'+1}{2}+\lfloor\frac{N-a'+2}{2}\rfloor .
   $$
If $N$ is even (resp. odd), then $\lfloor(N-a'+2)/2\rfloor=(N-a'+1)/2$ (resp.
$=(N-a'+2)/2$) and the claimed formula for $R(N)$ follows.
Now assume that $R(N)=\lfloor (N+1)/2\rfloor+1$ for some integer $N\in [1,2g-4]$. Let $i\in\N$ be such that $\gamma_i\leq N+1<\gamma_{i+1}$.
{}From Lemmas \ref{lemma3.4}(2) and \ref{lemma2.3}(2) we  have $\lfloor
(N+1)/2\rfloor+1=R(N)\leq i\leq (\gamma_i+2)/2\leq (N+3)/2$. If $N$ is
odd, then $i=(N+3)/2$ with $2\leq i\leq g-1$ and $\gamma_i=2i-2$, and (1)
follows from Remark \ref{rem2.2}. If $N$ is even, then $i=(N+2)/2$ with
$2\leq i\leq g-1$ and $\gamma_i\in\{2i-2,2i-1\}$, and (2) follows again
by Remark \ref{rem2.2}.
   \end{proof}
According to Remark \ref{rem3.1}, it is useful to compute $R(2g-2)$. The
result is the following.
   \begin{proposition}\label{prop3.3} Let $GS(\cX)=(\gamma_i:i\in\N)$ be
the gonality sequence of $\cX$. Then $R(2g-2)=\min\{2R(g-1),a\},$ where
   $$
a:=g-\max\{\gamma_i-(2i-2):i=1,\ldots,g\}\, .
   $$
   \end{proposition}
   \begin{proof} By Proposition \ref{prop3.1}, 
   $$
R(2g-2)=\min_{(i=1,\ldots,g)}\{2R(g-1),\ltilde(\gamma_i-1)+
\ltilde(2g-1-\gamma_i)\}\, .
   $$
   \begin{claim*} For $i=1,\ldots,g,$ $\ltilde(\gamma_i-1)=i-1$ and
$\ltilde(2g-1-\gamma_i)=g-\gamma_i+i-1.$
   \end{claim*}
In fact, the first claim follows immediately by the definition of the
function $\ltilde$. The proof of the second claim will follow from the
fact that $\gamma_{g-\gamma_i+i-1}<2g-1-\gamma_i<\gamma_{g-\gamma_i+i}\,
(*)$ for $i=1,\ldots,g$. To prove $(*)$, we apply induction on $i$. If
$i=1$, then $(*)$ becomes $\gamma_g=2g-2<2g-1<\gamma_{g+1}=2g$ by Lemma
\ref{lemma2.3}. Suppose that $(*)$ is true for $1\leq i<g$. We first show
that $2g-1-\gamma_{i+1}<\gamma_{g-\gamma_{i+1}+i+1}$. If this is not true,
then by inductive hypothesis and Proposition \ref{prop2.2} we would have
that 
   $$
\gamma_{g-\gamma_{i+1}+i+1}<2g-1-\gamma_{i+1}<2g-1-\gamma_i<
\gamma_{g-\gamma_i+i}\, .
   $$
Then $g-\gamma_{i+1}+i+1\leq g-\gamma_i+i-1$ so that
$\gamma_{i+1}-\gamma_i\geq 2$. Now let $a$ be an integer with
$\gamma_i<a<\gamma_{i+1}$. The $2g-1-a\in GS(\cX)$ by Proposition
\ref{prop2.2} and since $2g-1-\gamma_{i+1}<2g-1-a<2g-1-\gamma_i$ we must
have that $\gamma_{g-\gamma_{i+1}+i+1}<
\gamma_{(g-\gamma_i+i)-(\gamma_{i+1}-\gamma_i-1)}$, a contradiction. To
finish the proof of the claim we now prove the other inequality, namely 
$\gamma_{g-\gamma_{i+1}+i}<2g-1-\gamma_{i+1}$. If this were not true, then
we have
   $$
2g-1-\gamma_{i+1}<\gamma_{g-\gamma_{i+1}+i}<\gamma_{g-\gamma_{i+1}+i+1}\leq
\gamma_{g-\gamma+i}\, .
   $$
Then we must have $\gamma_{i+1}-\gamma_i\geq 2$, otherwise from the above
inequalities and inductive hypothesis it holds that
   $$
2g-1-\gamma_{i+1}=2g-2-\gamma_i
<\gamma_{g-\gamma_{i+1}+i}=\gamma_{g-\gamma_i+i-1}<
2g-1-\gamma_i\, ,
   $$
a contradiction. We now proceed as in the previous proof.
Now, since $\ltilde(\gamma_i-1)+\ltilde(2g-1-\gamma_i)=g-\gamma_i+2i-2$,
the result follows.
   \end{proof}
Now Theorem \ref{thm3.1}, Remark \ref{rem3.1} and the above computation of
$R(2g-2)$ imply the following.
   \begin{theorem}\label{thm3.2} Let $\cC=\cC(\cX,D,G)$ be an AG code such
that $2k\leq n$ and $n>2g,$ where $g$ is the genus of $\cX$. Let
$\gamma_2$ be the gonality of $\cX$ over $\fq.$ Then
    $$
s[\cC]\geq w(\cC)-g+\gamma_2-2\, .
    $$
   \end{theorem}
In the remaining part of this section we study the function $R$ on a
plane curve $\cX$ of degree $r+1$. In this case the genus of $\cX$ is
$g=r(r-1)/2$ and its gonality sequence $GS(\cX)$ is obtained from the
semigroup generated by $r$ and $r+1$ (cf. Lemma \ref{lemma3.3}(2)).
For an integer $a\in \N_0$, let $\alpha$ and $\beta$ be the non-negative
integers defined by 
    $$
a=\alpha r+\beta\, ,\qquad 0\leq \beta<r\, .
    $$
It is clear that $a\in GS(\cX)$ if and only $\beta\leq \alpha$.
    \begin{lemma}\label{lemma3.5} 
   $$
\ltilde(a)=\frac{\alpha(\alpha+1)}{2}+\min\{\alpha,\beta\}+1\, .
   $$
   \end{lemma}
   \begin{proof} If $a=0$, the formula is true so let $a>0$. Suppose first
that $a\in GS(\cX)$ so that $\min\{\alpha,\beta\}=\beta$. Then
$\ltilde(a)=1+2+\ldots+\alpha+\beta+1$ and we the claimed formula follows. 
Now let $a$ be a gap of $\cX$ so that $\beta>\alpha$. We have
$\ltilde(a)=\ltilde(\alpha r+\alpha)$ and the result follows by applying
the above computation to $\alpha r+\alpha\in GS(\cX)$.
   \end{proof}
   \begin{lemma}\label{lemma3.6} Let $N\in \N'\cap[-1,2g-2]$ and $a=\alpha
r+\beta$ a gap of $\cX$ with $\alpha\geq 1$ such that $a\leq N/2.$ Then
    $$
\ltilde(a)+\ltilde(N-a)\leq \ltilde(a-r)+\ltilde(N-(a-r))\, .
    $$
   \end{lemma}
   \begin{proof} Set $a':=a-r$, that is, $a'=(\alpha-1)r+\beta$. Let $b:=N-a=\delta r+\epsilon$, with $0\leq\epsilon<r$ so that 
$b'=N-a'=(\delta+1) r+\epsilon$. From Lemma \ref{lemma3.5} we have 
   $$
\ltilde(a)-\ltilde(a')=\alpha+1\, ,\qquad\text{and}\qquad
\ltilde(b)-\ltilde(b')\leq -\delta-1\, .
   $$
Now the result follows since $a\leq b=N-a$ implies $\delta\geq \alpha$.
    \end{proof}
Thus Proposition \ref{prop3.1} for the case of a plane curve becomes as follows.
    \begin{proposition}\label{prop3.4} For $N\in \N_0\cap[0,2g-2]$ let
$\alpha$ and $\beta$ be the integers defined by $\lfloor N/2\rfloor=\alpha
r+\beta$ with $0\leq \beta<r.$ Assume that $\alpha\geq 1.$
    \begin{enumerate}
\item[\rm(1)] If $\lfloor N/2\rfloor$ is a gap of $\cX,$ then
   $$
R(N)=\min\{\ltilde(\lfloor\frac{N}{2}\rfloor)+\ltilde(\lceil
\frac{N}{2}\rceil), \ltilde(\alpha r-1)+\ltilde(N-\alpha r+1) \}\, .
   $$
\item[\rm(2)] If $\lfloor N/2\rfloor\in GS(\cX),$ then
   $$
R(N)=\ltilde(\alpha r-1)+\ltilde(N-\alpha r+1)\, .
   $$
    \end{enumerate}
    \end{proposition}
    \begin{proof} It follows from Proposition \ref{prop3.1} and Lemma
\ref{lemma3.6}.
   \end{proof}
To improve this result we shall introduce the notion of ``jump". An
integer $N$ with $0\leq N\leq 2g-2$ is called {\em a jump of $\cX$}
whenever $R(N)>R(N-1)$ (so, $R(N)=R(N-1)+1$ by Lemma
\ref{lemma3.4}(1)). We denote by $U(\cX)$ the set of jumps of
$\cX$. Clearly $|U(\cX)|=R(2g-2)$ and this number can be computed via
the above proposition. More precisely the following holds.
    \begin{lemma}\label{lemma3.7} Let $\cX$ be a (non-singular) plane
curve of degree $r+1.$ Then
    \begin{enumerate}
\item[\rm(1)] $|U(\cX)|=\begin{cases}
r^2/4      & \text{if $r$ is even$,$}\\
(r^2-1)/4  & \text{if $r$ is odd$;$}
    \end{cases}$
\item[\rm(2)] $U(\cX)=\{\alpha r+\beta$:$ 
\text{$-1\leq \alpha\leq r-1,$ $0\leq \beta\leq r-1,$ and $2\beta+2\leq
\alpha$ or $\beta=r-1$}\}\setminus\{2g-1\}.$
    \end{enumerate}
    \end{lemma}
    \begin{proof} (1) Let us compute $R(2g-2)$. If  $r$ is even, then
$g-1=(r-2)(r+1)/2=(r-2)r/2+(r-2)/2$
and thus it belongs to $GS(\cX)$. By Proposition \ref{prop3.4},
$R(2g-2)=\ltilde(\alpha r+1)+\ltilde(2g-2-\alpha r+1)$ with
$\alpha:=(r-2)/2$. Now the result follows by applying Lemma
\ref{lemma3.5}. The case $r$ odd is similar.
(2) Let us denote by $T$ the set of the right-hand side in the equality in
Item (2). We claim that $|T|=R(2g-2)$. Indeed 
$|T|=\sum_{\beta=0}^{\lfloor(r-4)/2\rfloor}(r-2\beta-3)+r-1=R(2g-2)$.
Therefore it is enough to show that $T\subseteq U(\cX)$. {}From
Proposition \ref{prop3.4} and Lemma \ref{lemma3.5} it is easily seen that all
elements in $U(\cX)$ are jumps. Then the proof is complete.
    \end{proof}
Graphically, the set $U(\cX)$ looks like in the following example.
   \begin{example}\label{ex3.2} Let $\cX$ be a plane curve of degree $8$.
So $r=7$ and $g=21$. The next table shows all integers from $-1$ to
$2g-2=40$. The jumps of $\cX$ are marked in bold face.
    \begin{center}
    \begin{tabular}{ccccccc}
{\bf -1} & 0        & 1        & 2  & 3  & 4  & 5  \\
{\bf 6}  & 7        & 8        & 9  & 10 & 11 & 12 \\
{\bf 13} & {\bf 14} & 15       & 16 & 17 & 18 & 19 \\
{\bf 20} & {\bf 21} & 22       & 23 & 24 & 25 & 26 \\
{\bf 27} & {\bf 28} & {\bf 29} & 30 & 31 & 32 & 33 \\
{\bf 34} & {\bf 35} & {\bf 36} & 37 & 38 & 39 & 40
    \end{tabular}
    \end{center}
    \end{example}
Finally the promised improved description of $R(N)$ in the case of plane
curves is as follows.
    \begin{proposition}\label{prop3.5} Let $\cX$ be a (non-singular) plane
curve of degree $r+1$ and $N\in \N'\cap[-1,2g-2].$ Let $\alpha$ and $\beta$
be the integers defined by $N=\alpha r+\beta$ with $0\leq \alpha\leq r-2$
and $-1\leq \beta\leq r-2.$
    \begin{enumerate}
\item[\rm(1)] If $\beta>\lfloor\alpha/2\rfloor-1,$ then $R(N)=R(\alpha
r+\lfloor\alpha/2\rfloor-1);$
\item[\rm(2)] If $\beta\leq \lfloor\alpha/2\rfloor-1,$ then
     $$
R(N)=\begin{cases}
\alpha(\alpha+2)/4+\beta+1  & \text{if $\alpha$ is even,}\\
(\alpha+1)^2/4+\beta+1  &  \text{if $\alpha$ is odd.}
   \end{cases}
    $$
    \end{enumerate}
   \end{proposition} 
    \begin{proof} (1) In this case $\alpha r+\lfloor\alpha/2\rfloor-1$ is
the largest jump of $\cX$ not exceeding $N$ and (1) follows.
(2) Here place all the integers from $-1$ to $2g-2$ in an array according
to the corresponding values of $\alpha$ and $\beta$ (cf. Example
\ref{ex3.2}). The $j$-th row of the array contains $\lfloor
(j+2)/2\rfloor$ jumps of
$\cX$ which are precisely the ones in the first $\lfloor(j+2)/2\rfloor$
columns of the array. Thus the number of jumps from $-1$ to $N$ is:
$\beta+1$ in the row $\alpha$ plus
$2(\sum_{i=1}^{(\alpha-1)/2}i)=\alpha(\alpha+2)/2$ if $\alpha$ is even, and $\beta+1$ in  row $\alpha$ plus 
$2(\sum_{i=1}^{(\alpha-2)/2}i)=(\alpha+1)^2/4$ if $\alpha$ is odd. 
    \end{proof}
 
    \section{A property of Self-orthogonal codes}\label{s5}
In this section we state a new property of self-orthogonal codes. This property was first noticed by  Blackmore and Norton in \cite{blackmore-norton2} for the case of Hermitian codes, and used to improve their bounds on $s(\cC)$. Here we shall show that it holds in the very  general context of formally self-orthogonal codes.
We begin with some definitions. Let $\cC$, $\cC'$ be two codes of the same
length $n$ over $\fq$. We say that
they are {\em formally equivalent} (denoted $\cC\sim\cC'$) if there exists
an $n$-tuple $\bx$ of nonzero elements in $\fq$ such that $\cC=\bx*\cC$,
where $*$ stands for the coordinate-wise multiplication, see
\cite{munuera-pel}. The code $\cC$ is called {\em self-orthogonal} if
$\cC\subseteq \cC^{\perp}$, and {\em formally self-orthogonal} if there
exists an $n$-tuple $\bx$ of nonzero elements in $\fq$ such that
$\cC\subseteq\bx*\cC^{\perp}$.  
    \begin{example}\label{ex51} (1) Let
$\cC=RM_q(r,m)$ be a $q$-ary Reed-Muller code. Since
$RM_q(r,m)^{\perp}=RM_q(m(q-1)-r-1,m)$, then $\cC$ is self-orthogonal
whenever its dimension is at most $n/2$. 
(2) Let $\cC=\cC(\cX,D,G)$, $\cC'=\cC(\cX,D,G')$ be two AG code as the
ones treated
in the former sections. Then $\cC\sim\cC'$ if and only if $G\sim G'$,
where $\sim$ stands for the usual equivalence of divisors, see
\cite{munuera-pel}. The dual of $\cC$ is $\cC^{\perp}=\cC(\cX,D,D+W-G)$,
where $W$ is a canonical divisor obtained as the divisor of a differential
form having simple poles and residue 1 at every point in $\supp(D)$. Thus
we deduce that $\cC$ is formally self-orthogonal if there is an effective
divisor $E$ such that $\supp(E)\cap\supp(D)=\emptyset$ and $D+W-2G\sim E$.  
Now, let $\cC=\cC(\cX,D,G)$ be a Hermitian code, that is, a code constructed from a Hermitian curve  of affine equation $y^q+y=x^{q+1}$ over ${\mathbf F_{q^2}}$, by taking $Q:=(0:1:0)$, $D$ equals to the sum of the $q^3$ affine points, and $G=mQ$ (see \cite{blackmore-norton2}, \cite{yang-kumar-sti}). Since
$D+W-2G\sim (n+2g-2-2m)Q$, then $\cC$ is self-orthogonal whenever its
dimension is at most $n/2$. 
    \end{example}
We shall give a bound on the trellis state complexity of formally
self-orthogonal codes. This bound is based on the following result. 
    \begin{proposition}\label{fs-o} Let $\cC$ be a formally
self-orthogonal $[n,k,d]$ code$.$ Then for $1\le i\le n$ either
$\cP_{i-1}=
\cP_{i}$ or $\cF_{i-1}= \cF_{i}.$
    \end{proposition} 
    \begin{proof} Let $\bx$ be such that $\cC\subseteq\bx*\cC^{\perp}$. If
$\cP_i$ and $\cF_{i+1}$ were proper subspaces of $\cP_{i+1}$ and
$\cF_i$ respectively, then there exist codewords
$c=(c_1,\dots,c_i,0\cdots,0)$, $c'=(0,\dots,0,c'_i,\dots,c'_n)$ in
$\cC$
such that $c_ic'_i\neq 0$. On the other hand, being $\cC$ self-orthogonal
we have $0=(\bx*c)\cdot c'=x_ic_ic'_i$, hence $c_ic'_i=0$, which is a
contradiction. 
   \end{proof}
   \begin{example}\label{ex5.2} Let us see how the
above result can be used to improve the bounds on $s[\cC]$. Let $\cC$ be a
$[n,k,d]$ code with $2k\le n$. As we know we have $s[\cC]=k$ if $2d\ge
n+2$ (Proposition \ref{prop2.1}). Set $p_i:=\dim_\fq\cP_i$ and
$f_i:=\dim_\fq\cF_i$. Let us examine the border case $2d=n+1$.
Again according to Proposition \ref{prop2.1}, $\Delta(\cC)=\min\{
p_{d-1}+f_{d-1},
p_{d}+f_{d} \}$ hence $\Delta(\cC)\le 1$ as $p_{d-1}=f_d=0$,
$p_{d},f_{d-1}\le 1$. If $\cC$ is formally self-orthogonal then either
$p_{d}=0$ or $f_{d-1}=0$ and thus $s[\cC]=k$.
   \end{example}
   \begin{proposition}\label{cor5.1} Let $\cC$ be a formally self-orthogonal
$[n,k,d]$ code with $2d\leq n+1.$ Then for $0\le i\le j\le n$ we have
$(p_{j}-p_i)+(f_i-f_{j})\le
j-i.$ As a consequence$,$ $\Delta[\cC]\le \lfloor\frac{n-2d+2}{2}\rfloor.$ 
    \end{proposition} 
    \begin{proof} Apply Proposition \ref{fs-o} $j-i$ times. By
this result with $i=d-1$ and $j=n-d+1$, we obtain $f_{d-1}+p_{n-d+1}\le
n-2d+2$, hence either $\Delta_{d-1}\le \lfloor(n-2d+2)/2\rfloor$ or
$\Delta_{n-d+1}\le \lfloor(n-2d+2)/2\rfloor$. 
    \end{proof} 
   
    \end{document}